\documentclass[a4paper]{amsart}
\title{The Witten top Chern class via $K$-theory}

\author{Alessandro Chiodo}
\address{Alessandro Chiodo\\ Laboratoire J.-A.~Dieudonn\'e\\ Universit\'e de Nice\\ Parc \linebreak Valrose\\ 06108 Nice cedex 2\\ France}
\email{chiodo@math.unice.fr}
\thanks{Supported by
the Istituto Nazionale di Alta Matematica,
and the Marie Curie Intra-European
Fellowship within the 6th European Community Framework Programme,
MEIF-CT-2003-501940.}

\normalsize \linespacing=\baselineskip
\date{}

\dedicatory{}

\urladdr{http://math.unice.fr/$\sim$chiodo/}

\usepackage{amsmath, amsthm, amssymb, setspace, graphics}
\usepackage{mathrsfs}%for some cursive fonts
\usepackage[all]{xy}

\numberwithin{equation}{subsection}
\swapnumbers
\newtheorem{thm}[equation]{Theorem}

\newtheorem{lem}[equation]{Lemma}
\newtheorem{pro}[equation]{Proposition}
\theoremstyle{remark}
\newtheorem{defn}[equation]{Definition}
\newtheorem{rem}[equation]{Remark}

\newcommand{\CC}{\mathbb C}
\newcommand{\EE}{\mathbb E}

\newcommand{\PP}{\mathbb P}
\newcommand{\QQ}{\mathbb Q}

\newcommand{\VV}{\mathbb V}
\newcommand{\ZZ}{\mathbb Z}

\newcommand{\id}{\operatorname{id}}

\newcommand{\im}{\operatorname{im}}

\newcommand{\codim}{\operatorname{codim}}

\newcommand{\Hom}{\operatorname{Hom}}

\newcommand{\Spec}{\operatorname{Spec}}

\newcommand{\Vect}{\operatorname{Vect}}
\newcommand{\al}{\alpha}

\newcommand{\fie}{\varphi}
\newcommand{\ka}{\kappa}
\newcommand{\la}{\lambda}

\newcommand{\ctop}{c_{\mathrm{top}}}
\newcommand{\rk}{\operatorname{rk}}

\newcommand{\class}{\operatorname{cl}}
\newcommand{\ke}{\operatorname{Ke}}

\newcommand{\ch}{\operatorname{ch}}
\newcommand{\td}{\operatorname{td}}

\newcommand{\Sym}{{\operatorname{Sym}}}
\newcommand{\calCh}{\operatorname{\bf{Ch}}}
\DeclareMathAlphabet{\mathpzc}{OT1}{pzc}{m}{it}
\newcommand{\stack}{\mathscr}

\newcommand{\typek}{{\pmb k}}
\newcommand{\nmarkings}{s_1,\dots,s_n}

\newcommand{\stC}{\stack C}
\newcommand{\stL}{\stack L}

\def\ol{\overline}

\def\wt{\widetilde}

\newcommand{\SSS}{\stack S}

\newcommand{\MM}{\stack M}

\newcommand{\MS}{{{\SSS}}_{g,n}({r,\typek})}

\newcommand{\MSbar}{\ol{{\SSS}}_{g,n}({r,\typek})}

\newcommand{\cvirt}{c_W}
\newcommand{\kb}{{\pmb k}}

\begin{document}

\begin{abstract}
{The Witten top Chern class is the crucial
cohomology class needed to state a conjecture
by Witten relating
the Gelfand--Diki\u\i\  hierarchies
to higher spin curves.
In \cite{PV},
Polishchuk and Vaintrob provide an algebraic construction
of such a class.
We present a more straightforward
construction via $K$-theory.
In this way we
short-circuit the passage through
bivariant intersection theory and the use of
MacPherson's graph construction. Furthermore,
we show that the Witten top Chern class
admits a natural lifting
to the $K$-theory ring.
}
\end{abstract}

\maketitle

\pagestyle{myheadings}

\markboth{ALESSANDRO CHIODO}{The Witten top Chern class via $K$-theory}

\section{Introduction}
\subsection{Witten's conjecture.}
In \cite{Wir}, Witten conjectures that
the intersection numbers of certain cohomology classes
on the moduli stacks of stable $r$-spin curves
encode a solution of the Gelfand--Diki\u\i\
(also known as the higher KdV) hierarchy. (The conjecture
is a generalization of the Kontsevich--Witten Theorem
\cite{Wi2}, \cite{Ko}.)
Witten's formulation of the conjecture lacks both a
rigorous definition of the moduli stack $\MSbar$ of
stable  $r$-spin curves
and a construction of the
crucial cohomology class on $\MSbar$
which is usually referred
to as the \emph{Witten top Chern class}.

\subsection{Definition of the moduli stack.}
The moduli stack of
stable $r$-spin curves was first defined by Jarvis \cite{Jageom}
and is the compactification of the stack $\MS$
labeled by the integer and  nonnegative
indexes $r, g,n,$ and $\kb=(k_1,k_2,\dots,k_n)$
satisfying
$r\ge 2$, $2g-2+n>0$, and $2g-2-\sum_i k_i\in r\ZZ$.
The stack $\MS$ classifies smooth stable $r$-spin curves:
$n$-pointed smooth stable curves $(C; \nmarkings)$ of genus $g$
equipped with a line bundle
$L$ on $C$ and an isomorphism $f\colon
L^{\otimes r} \to \omega_C(-\sum_i k_i [s_i])$.
The compactified stack $\MSbar$ classifies
stable curves $C$ equipped
with a torsion free sheaf of rank one $L$
and a nonzero homomorphism
$f$ as above. (Using stack-theoretic curves,
this solution
is rephrased in \cite{AJ} and
can be modified as
shown in \cite{Chmod};
however, the problem of constructing
the Witten top Chern class
has equivalent solutions with all
these compactifications, Remark \ref{compstack}.)

\subsection{Witten's definition of his
class on an open substack.}
In \cite{Wir},
the Witten top Chern class is defined
on the open substack $\stack V$ of $\MS$
satisfying
\begin{align}\label{opensub}
H^0(C_x,L_x)=0 \quad \text{for all}\quad x\colon\Spec \CC\to \stack V,
\end{align}
where $((C_x; s_{1,x},\dots,s_{n,x}),L_x,f_x)$ is the
$r$-spin curve represented by $x$.
Over $\stack V$ there exists a vector bundle whose fibre
on the point $x$ is the space \linebreak $H^1(C_x,L_x)$.
This is the higher direct image $R^1\pi_*\stack L$,
where $\pi\colon \stack C\to \stack V$
is the universal curve
and $\stack L$ on $\stack V$
is the universal line bundle.
Then we take
\begin{equation}\label{wittensdef}
\ctop (R^1\pi_*\stack L).
\end{equation}
If $g=0$, the condition \eqref{opensub} is satisfied over the
whole moduli stack.
In general, $R^1\pi_*\stL$ is not a vector bundle,
so the usual definition of top Chern class
does not apply.
Witten sketches a generalization based on the index
theory of elliptic operators. It is not clear, however,
how to extend this approach to singular curves and
to the whole stack $\MSbar$.
\subsection{Polishchuk and Vaintrob's construction.}
In \cite{PV}, Polishchuk and Vaintrob
provide an algebraic construction of the Witten
top Chern class on the whole stack $\MSbar$.
They work with the universal
torsion free sheaf of rank one
$\stL$ on the universal family $\stC\to \MSbar$, and
they consider the pushforward
$R\pi_*\stL$ in the derived category, which is
used to construct a $\ZZ/2\ZZ$-graded
spinor bundle $S$. By extending
MacPherson's graph construction
to $2$-periodic complexes, Polishchuk and Vaintrob define
the localized Chern character
of $S$ in bivariant intersection theory.
Passing to the Chow ring,
they obtain a Chow cohomology class $c_{PV}$
(see Section \ref{par:pv}).
Let $\chi$ denote $\chi(\stC,\stL)$; we have
$c_{PV}\in A^{-\chi}(X)_{\QQ}$.
Furthermore, such
a class is compatible with Witten's earlier definition and,
by \cite{Po}, satisfies all the axioms of
the cohomological field theory defined in \cite{JKV}---this
is a preliminary condition
to Witten's conjecture \cite{Wir}.

\subsection{Our construction of the Witten top Chern class.}
We start from the pushforward
$R\pi_*\stL$ in the derived category, but follow a
different path, which is more straightforward and
closer to Witten's original idea \eqref{wittensdef}
(our approach is also alluded to in
\cite[\S2, p.2, Remark]{Po}).
We obtain the class $c_{PV}$ by working only with classes
in $K$-theory, short-circuiting the passage through
$\ZZ/2\ZZ$-graded spinor bundles,
bivariant intersection theory,
and the use of MacPherson's
graph construction.

We define the \emph{$K$-theory Euler class}, a
generalization of the
total lambda class evaluated at $-1$,
which by definition is
\begin{align*}
\lambda_{-1}\colon \Vect(X)&\to K_0(X) \\
V &\mapsto \sum\nolimits_i(-1)^i[\Lambda^i{V}].
\end{align*}
The $K$-theory Euler class $\ke(F_\bullet,a)$
is defined for a pair
$(F_\bullet,a)$ where $F_\bullet$
is a bounded complex in degrees $1$ and $0$
of locally free coherent sheaves and 
$a\colon \mathcal O\to \Sym^{r-1} F_0\otimes F_1$ is
a closed and nondegenerate form
(see Definition
\ref{closed}
and Definition \ref{nondeg}).

Once $\ke(F_\bullet,a)$ is defined,
then, the construction
follows naturally.
We apply it to a double complex $F_\bullet$ in
degrees $1$ and $0$ representing  $(R\pi_*\stL)^\vee$
in the derived category of $\MSbar$.
Indeed, $F_\bullet$ is equipped with
a closed and nondegenerate form $a$
defined using $\stL^{\otimes r}\to \omega$ (Section \ref{sect:wtcc}).

For a vector bundle $V$, the total lambda class
$\la_{-1}$ is related
to the top Chern class in the
Chow ring by
\begin{equation}\label{ke-ctop}
c_{\rm top}(V)= \ch(\la_{-1}(V^\vee))\cdot \td (V).\end{equation}
In Definition \ref{defnctop} and
Definition \ref{wtcc-chow}, we
define a Chow cohomology class $c_W$ by
\begin{equation}
c_W=\ch(\ke(F_\bullet,a))\frac{\td (F_1^\vee)}{\td({F_0^\vee})}.\end{equation}
In Theorem \ref{thm:compPV},
we prove the identity $c_W=c_{PV}$.

It is also worth mentioning that
our construction
produces
a $K$-class $K_W$ that lifts $c_W$ to the
$K$-theory ring, Definition \ref{wtcc-chow}.
We hope that such a lifting can
not only clarify the definition of the Witten top Chern class,
but also improve our understanding of the conjecture.

\subsection{Structure of the paper.}
In Section \ref{sect:notn}, we introduce our notation.
In Section \ref{sect:euler}, we define the $K$-theory Euler class.
We provide our construction for the Witten top
Chern class in Section \ref{sect:wtcc} and
prove the identity with the Polishchuk--Vaintrob class in
Section \ref{sect:pv}.
Finally, in Section \ref{sect:comp}, we
give some explicit examples.

\subsection{Acknowledgments.}
I would like to thank Constantin Teleman, who was most generous to
me; I had several illuminating discussions with him that put this
project on the right track.
I am also extremely grateful to my Ph.D. supervisor
Alessio Corti for shaping my views on this problem and
for putting a great deal of care
in assisting in the preparation of
this paper.
At the Laboratoire Dieudonn\'e in Nice
I benefited from many discussions with the members of
the \'Equipe G\'eom\'etrie Alg\'ebrique.
In particular I would like to
thank Charles Walter for
several explanations and for
pointing out Lemma \ref{lem:K-2-per}.
\section{Notation}\label{sect:notn}
\subsection{Schemes.}
All schemes are of finite type over $\CC$.

\subsection{Fibres.}
\label{fibres:notn}
Let $p\colon E\to X$ be a
complex vector bundle on a scheme;
for $x\in X$ we denote by $E_x$ the fibre
$p^{-1}(x)$ over $x\in X$.
For a morphism of vector bundles $f\colon E\to F$,
we denote by $f_x\colon E_x\to F_x$ the morphism
induced on the fibre over $x\in X$.
For coherent locally free sheaves $K$ on $X$ we
denote by $K_x$ the fibre of the corresponding
vector bundle.

\subsection{$K$-theory.}
We denote by $\Vect(X)$ the category of finite
complex vector bundles on $X$.
We denote by $K_0(X)$ and $K'_0(X)$
the Grothendieck
groups generated by coherent locally free sheaves
and coherent sheaves on $X$.

\subsection{Symmetric and exterior product homomorphisms.}
Let $E$ be a coherent locally free sheaf on a scheme $X$.
For any nonnegative integers $h$ and $k$, the
homomorphisms
\begin{equation}\label{products}
\sigma_n\colon {\Sym}^{h}{E}\otimes{\Sym}^{n} {E}\rightarrow
\Sym^{h+n}{E}\quad{\rm and} \quad
\lambda_n\colon\Lambda^h{E}\otimes\Lambda^{n}{E} \rightarrow
\Lambda^{h+n} {E}
\end{equation}
are the natural {\em symmetric product} and {\em exterior
product}.

\subsection{Symmetric and exterior powers of a complex.}
\label{symantisym}
In this paper we work with complexes in degree $0$
and $1$; in this case, the symmetric and the exterior power
in the derived category can be realized by bounded
cohomological and homological Koszul complexes.
Following \cite[\S1]{GrKos}, we
write the complexes explicitly.

First, for any coherent locally free sheaf $E$
on $X$ and for any nonnegative integer $h$,
denote by
$$s\colon \Sym^h E\to \Sym^{h-1}E\otimes E
\quad \text{and} \quad l\colon \Lambda^h E\to \Lambda^{h-1}E\otimes E,$$
the dual homomorphisms of the products
$\sigma_1$ and $\lambda_1$ applied to
$E^\vee$. This allows us to define, for any
$d\colon V^0\to V^1$,
the homomorphisms
\begin{multline}\label{eq:kgoesup}\Sym^h V^0 \otimes \Lambda^k {V^1} \to
\Sym ^{h-1}V^0 \otimes V^0 \otimes \Lambda^{k}{V^1}\\
\to \Sym ^{h-1}V^0 \otimes V^1 \otimes \Lambda^{k}{V^1}\to
\Sym^{h-1}{V^0}\otimes \Lambda^{k+1}{V^1},\end{multline}
where the first homomorphism
is $s$ tensored with the identity on $\Lambda ^k{V^1}$,
the second homomorphism is the identity on $\Sym^{h-1}{V^0}$
and $\Lambda ^{k}{V^1}$ tensored with $d$,
and
the third homomorphism is the
identity on $\Sym^{h-1}{V^0}$ tensored with the exterior product
$\lambda_1$.
In this way, for any complex of
coherent locally free sheaves
$V^\bullet: 0\to V^0\to V^1\to 0$ we get
the {\em symmetric
power complex} $\Sym^N(V^\bullet)$
\begin{multline}
\label{sym} 0\rightarrow \Sym^N V^0\rightarrow
\Sym^{N-1}V^0\otimes V^1\rightarrow
\Sym^{N-2}V^0\otimes \Lambda^2 V^1\rightarrow \cdots \\
\cdots\rightarrow \Sym^{2}V^0\otimes \Lambda^{N-2} V^1 \rightarrow
V^0\otimes \Lambda^{N-1} V^1 \rightarrow \Lambda^N V^1\rightarrow 0.
\end{multline}

Similarly,
for any complex
$W_\bullet:0\to W_1\xrightarrow{} W_0\to 0$
of coherent locally free sheaves on $X$.
We can define
\begin{multline}\label{eq:hgoesup}\Sym^h {W_0} \otimes
\Lambda^k {W_1} \to
\Sym ^{h}{W_0} \otimes {W_1} \otimes \Lambda^{k-1}{W_1}\\
\to \Sym ^{h}{W_0} \otimes {W_0} \otimes \Lambda^{k-1}{W_1}\to
\Sym^{h+1}{W_0}\otimes \Lambda^{k-1}{W_1}\end{multline}
as the composite
$(\id\otimes l)\circ(\id\otimes
d\otimes\id)\circ({\sigma_1\otimes \id})$.
The {\em exterior
power complex} $\Lambda^N(W_\bullet)$
is the complex
\begin{multline}
\label{antisym} 0\rightarrow \Lambda^N W_1\rightarrow
W_0\otimes\Lambda^{N-1}W_1\rightarrow
\Sym^{2}W_0\otimes \Lambda^{N-2} W_1\rightarrow \cdots \\
\cdots\rightarrow \Sym^{N-2}W_0\otimes \Lambda^{2} W_1
\rightarrow\Sym^{N-1}W_0\otimes W_1 \rightarrow \Sym^N
W_0\rightarrow 0.
\end{multline}

In fact, in characteristic
$0$, $\Sym^N$ and $\Lambda^N$
are well defined functors up to
homotopy \cite{De}. We denote by $\Sym^N(f)$ and
$\Lambda^N(f)$ the homomorphism of complexes associated to a
homomorphism $f$ of complexes in degrees $1$ and $0$.

\begin{pro}\label{cohoantisym}
Let $d\colon W_1\to W_0$ be a homomorphism
of complex vector spaces, and let $H_1$
and $H_0$ be its the kernel and its
cokernel. The cohomology of
$\Lambda^{N}(W_\bullet)$ is given by
$\Sym^{h}H_0\otimes \Lambda^{N-h}H_1$ for $0\le h\le N$.
\end{pro}
\begin{proof}
Let $V$ be the image of $d$.
The complex
$W_\bullet$ is the sum
of $H_1\xrightarrow {0}H_0$ and
$V\xrightarrow{\id} V$.
Therefore, we get
$$\Lambda^N(W_\bullet)=\sum_{i=0}^N
\Lambda^{N-i}(H_\bullet) \otimes \Lambda^i(V\xrightarrow{\id}
V).$$
Note that
$\Lambda^i(V\xrightarrow{\id} V)$
is exact for $i>0$.
The claim follows.
\end{proof}

\section{The $K$-theory Euler class}\label{sect:euler}
The total lambda class evaluated at $-1$ is defined
over vector bundles up to isomorphism:
\begin{align*}
\lambda_{-1}\colon \Vect(X)&\to K_0(X) \\
V &\mapsto \sum\nolimits_i(-1)^i[\Lambda^i{V}].
\end{align*}
We now extend its definition.
We choose an integer
$$m\ge 1.$$
We work with complexes
$F_\bullet=(0\to F_1 \to  F_0\to 0)$
of coherent locally free sheaves $F_i$ on the scheme $X$.
Denote by $H_{1,x}$ and $H_{0,x}$ the kernel and the
cokernel of $(F_\bullet)_x$.

\subsection{Closed and nondegenerate forms.}
Consider a homomorphism
$$a\colon {\mathcal O}_X\to \Sym^m{F_0}\otimes{F_1}.$$
\begin{defn}\label{closed}
The form $a$ is {\em closed}
if $d\circ a=0$ where $d$ is the natural map
$d\colon\Sym^m{F_0}\otimes{F_1}\to \Sym^{m+1}F_0$.
\end{defn}
\begin{rem}\label{expl-closed}
Consider the complex
$\Lambda^{m+1}(F_\bullet)$
\begin{multline}\label{antisymWTCC}
0\xrightarrow{d} \Lambda^{m+1} {F_1}\xrightarrow{d}
{F_0}\otimes\Lambda^{m}{F_1}\xrightarrow{d}
\Sym^{2}{F_0}\otimes \Lambda^{m-1} {F_1}\xrightarrow{d} \cdots \\
\cdots\xrightarrow{d} \Sym^{m-1}{F_0}\otimes \Lambda^{2}
{F_1} \xrightarrow{d}\Sym^{m}{F_0}\otimes {F_1}
\xrightarrow{d} \Sym^{m+1} {F_0}\xrightarrow{d} 0,
\end{multline}
introduced in \eqref{antisym}. Note that the form $a$ is closed if and
only if $d(a(1))$ vanishes.

By Proposition
\ref{cohoantisym}, the cohomology
of the complex $\Lambda^{m+1}(F_\bullet)_x$
at \linebreak $(\Sym^{m}{F_0}\otimes {F_1})_x$ is
$\Sym^{m}H_{0,x}\otimes H_{1,x}$.
For
any $x\in X$, if $a$ is closed,
$a(1)$ induces an element in
$\Sym^{m}H_{0,x}\otimes H_{1,x}$, which can be regarded
as a  linear system
\begin{equation}\label{linsyst}
S_x\colon H_{1,x}^\vee\rightarrow \Sym^m H_{0,x}=
H^0(Y,{\mathcal O}_Y(m))
\end{equation}
with $Y=\PP H_{0,x}^\vee$.
\end{rem}
\begin{defn}\label{nondeg}
A closed form $a$ is {\em nondegenerate} if for any $x\in X$ the
image of $S_x$ is a base point free linear system on $\PP
H_{0,x}^\vee$.
\end{defn}
\begin{rem}\label{expl-nondeg1}
Note that this condition implies
$\rk(H_{0,x}) \le \rk(H_{1,x})$
and, therefore,  $\rk(F_0) \le \rk(F_1)$.
\end{rem}
\begin{rem}\label{expl-nondeg2}
Note that a closed form $a$ is nondegenerate if and only if for
any $x\in X$ the following condition is satisfied:
an element $v\in H_{0,x}^\vee$ is zero if
\begin{equation}\label{condnondeg}
\forall w\in H_{1,x}^\vee, \quad \langle S_x(w),v^m\rangle=0.
\end{equation}

\end{rem}

\subsection{The double complex $L^{\bullet,\bullet}$.}
From now on we always assume that $a$ is a closed and
nondegenerate form. We write
\begin{equation}\label{complexL}L^{h,k}=\Sym^h {F_0}\otimes\Lambda^k {F_1}.
\end{equation}
We define homomorphisms
$$\wt d\colon L^{h,k}\to L^{h+1,k-1},$$
$$\wt a\colon L^{h,k}\to L^{h+m,k+1}$$
of bidegrees $(1,-1)$ and $(m,1)$. The homomorphism $\wt d$ is
defined as the differential of the
complex $\Lambda^N(F_\bullet)$
at \eqref{antisym} with $N=h+k$. Consider the natural
homomorphisms described in \eqref{products}
\begin{align*}
\sigma_m\colon&\Sym^{h} {F_0}\otimes \Sym^m {F_0}\rightarrow
\Sym^{h+m} {F_0}, & \lambda_1\colon& {F_1}\otimes {\Lambda}^k
{F_1} \rightarrow {\Lambda}^{k+1}{F_1}.
\end{align*}
The homomorphism $\wt a$ is the composite of ${\rm id}\otimes
{a}\otimes{\rm id}$ and $\sigma_m\otimes  \lambda_1$
\begin{multline}
\label{afactors}
\Sym^h {{F_0}}\otimes \Lambda^{k} {{F_1}}
\xrightarrow{\id\otimes a\otimes\id} \Sym^{h} {{F_0}}\otimes
\Sym^{m}{{F_0}}\otimes {{F_1}}\otimes \Lambda^{k} {{F_1}}\\
\xrightarrow{\sigma_m\otimes \lambda_1} \Sym^{h+m} {{F_0}}\otimes
\Lambda^{h+1} {{F_1}}.
\end{multline}

\begin{lem}\label{doublecx}
The bigraded sheaf $L^{\bullet,\bullet}$ is a double complex with
differentials $\wt d$ and $\wt a$.
\end{lem}
\begin{proof}
We show $\wt d\circ \wt d=\wt a\circ \wt a=\wt d \circ \wt a+
\wt a\circ\wt d=0$.
The homomorphism $\wt d$ is the
differential of the exterior power
of $F_\bullet$; therefore $\wt d\circ \wt d=0$.
Furthermore, a local  description of $\wt a$ at $x\in X$
shows $\wt a\circ \wt a=0$.
We choose an open $U\ni x$ such that
$$a(1)=\sum\nolimits_{l=0}^N P_l\otimes
\sigma_l\in\Sym^{m}F_0(U)\otimes F_1(U)$$
with $P_l\in \Sym^m F_0(U)$ and $\sigma_l\in F_1(U).$
Then, $\wt a_x\colon L^{h,k}_x\to L^{h+m,k+1}_x$ can be written, for $c_i\in
(F_0)_x$ and $b_j\in (F_1)_x$, as the homomorphism  sending
$\prod\nolimits_{1\le i\le h}c_i\otimes \bigwedge_{1\le j\le k} b_j\in
L_x^{h,k}$ to
\begin{equation}
\label{aonfibre} \sum\nolimits_{l=0}^N\prod\nolimits_{1\le i\le h}c_i\cdot
P_l(x)\otimes \sigma_l(x)\wedge\bigwedge\nolimits_{1\le j\le k} b_j.
\end{equation}
Using $\sigma_l\wedge
\sigma_{l'}+\sigma_{l'}\wedge \sigma_l=0$
we see that $\wt a$ is a
differential.

Finally, we show $\wt d \circ \wt a+
\wt a\circ\wt d=0$.
We point out that $\wt d_x\colon L^{h,k}_x\to L^{h+1,k-1}_x$
sends
$\prod\nolimits_{1\le i\le h} c_i\otimes 
\bigwedge\nolimits_{0\le j\le k} b_j$ to
\begin{equation}\label{donfibre}
\sum\nolimits_{j_0=1}^{k}(-1)^{j_0}\left (\prod\nolimits_{i}c_i\cdot
d_x\left(b_{j_0}\right)\otimes \bigwedge\nolimits_{j\ne j_0} b_j\right),
\end{equation}
with $d_x=(d_{F_\bullet})_x$, $c_i\in(F_0)_x$, and $b_j\in
(F_1)_x$.
Indeed, as illustrated in Section \ref{symantisym},
the homomorphism $d_x$ is the composite
\begin{multline}
\Sym^h(F_{0,x})\otimes \Lambda^{k}(F_{1,x})\to
\Sym^h(F_{0,x})\otimes F_{1,x}\otimes \Lambda^{k-1}(F_{1,x})
\to \\ \to \Sym^h(F_{0,x})\otimes F_{0,x}\otimes \Lambda^{k-1}(F_{1,x})
\to \Sym^{h+1}(F_{0,x})\otimes \Lambda^{k-1}(F_{1,x}),\end{multline}
where the first arrow denotes
the identity homomorphism tensored by \linebreak $l_x\colon
\Lambda^{k}(F_{1,x})\to
F_{1,x}\otimes \Lambda^{k-1}(F_{1,x})$,
which we write as
$$\bigwedge\nolimits_{0\le j\le k} b_j\mapsto
\sum\nolimits_{j_0=1}^{k}(-1)^{j_0}\left(
b_{j_0}\otimes \bigwedge\nolimits _{j\ne j_0} b_j\right).
$$
Then the definition given in \eqref{eq:hgoesup}
implies \eqref{donfibre}.
As a consequence of \eqref{donfibre} and \eqref{aonfibre},
the fact that $\wt d$ and $\wt a$ anticommute is
just another way to say that $a$ is closed.
\end{proof}
\subsection{The complexes $K_i^{\bullet,\bullet}$.}
\label{doubleK}
For each $i=0,\dots,m$, set
$$K_i^{p,q}=\Sym^{p+mq-i}{F_0}\otimes\Lambda^{-p+q}{F_1}$$
with differentials of bidegrees $(1,0)$ and $(0,1)$
\begin{align*}
\wt d\colon K_i^{p,q}\rightarrow K_i^{p+1,q}
\quad\quad  \text{and} \quad \quad
\wt a\colon K_i^{p,q}\rightarrow K_i^{p,q+1}.
\end{align*}
This diagram illustrates $K_0^{\bullet,\bullet}$.
\[ \xymatrix@C=0.4cm{
    &
      \cdots &
            \Sym^{3m+1}{F_0}\otimes\Lambda^{2}{F_1} &
                \cdots \\
    \Sym^{2m-1}{F_0}\otimes\Lambda^{3}{F_1}\ar[r]^{\widetilde d}&
        \Sym^{2m}{F_0}\otimes\Lambda^{2}{F_1}\ar[r]^{\widetilde d}&
            \Sym^{2m+1}{F_0}\otimes {F_1}\ar[r]\ar[u]^{\widetilde a}&
                \dots \\
%    \Sym^{3m-1}{F_0}\otimes\Lambda^{4}{F_1}\ar[r]^{\widetilde d}\ar[u]^{\widetilde a}&
%        \Sym^{3m}{F_0}\otimes\Lambda^{3}{F_1}\ar[r]^{\widetilde d}\ar[u]^{\widetilde a}&
%            \Sym^{3m+1}{F_0}\otimes \Lambda^{2}{F_1}\ar[r]\ar[u]^{\widetilde a}&
%                \dots \\
%    \Sym^{2m-1}{F_0}\otimes\Lambda^{3}{F_1}\ar[r]^{\widetilde d}\ar[u]^{\widetilde a}&
%        \Sym^{2m}{F_0}\otimes\Lambda^{2}{F_1}\ar[r]^{\widetilde d}\ar[u]^{\widetilde a}&
%            \Sym^{2m+1}{F_0}\otimes {F_1}\ar[r]\ar[u]^{\widetilde a}&
%                \dots \\
    \Sym^{m-1}{F_0}\otimes\Lambda^{2}{F_1}\ar[u]^{\widetilde a}\ar[r]^{\widetilde d}&
        \Sym^{m}{F_0}\otimes {F_1}\ar[r]^{\widetilde d}\ar[u]^{\widetilde a}&
            \Sym^{m+1}{F_0}\ar[r]\ar[u]^{\widetilde a}&
                0\ar[u] \\
    0 \ar[u]\ar[r]&
        {\mathcal O_X}\ar[u]^{\widetilde a} \ar[r] &
            0 \ar[u]\ar[r] & 0\ar[u]
                \\
%0\ar[u]\ar[r]&0\ar[r]\ar[u]&0\ar[u] \ar[r]& 0\ar[u]\\
} \]
Fig.~1 illustrates that $K_i^{\bullet,\bullet}$
is nonzero only if
$0\le -p+q\le \rk(F_1)$ and $p+mq\ge 0$.
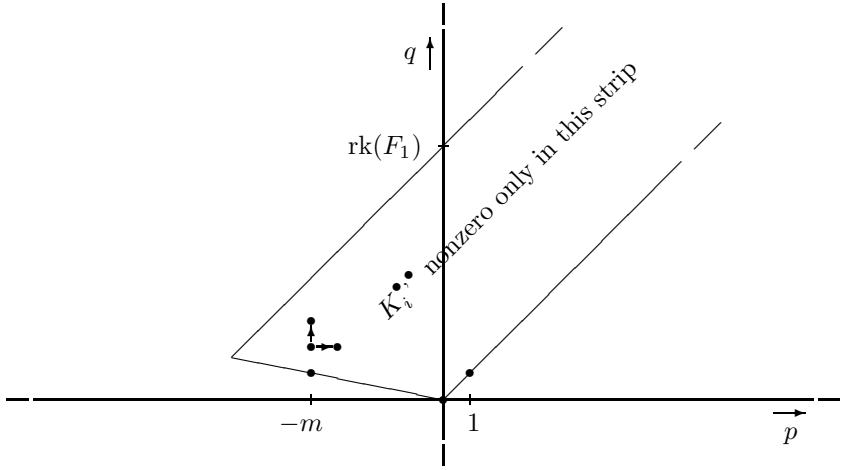
\begin{figure}[h]
\begin{picture}(250,200)(70,30)
\put(35,60){\line(1,0){8}} \put(45,60){\line(1,0){295}}
\put(342,60){\line(1,0){8}} \put(200,35){\line(0,1){8}}
\put(200,45){\line(0,1){155}} \put(200,202){\line(0,1){8}}
 \put(150,62){\line(0,-1){4}}
 \put(210,62){\line(0,-1){4}}
 \put(209,48){$1$}
 \put(138,48){$-m$}
 \put(185,189){$q$}
 \put(195,185){\vector(0,1){12}}
 \put(329,45){$p$}
 \put(325,55){\vector(1,0){12}}
 \put(164,153){$\rk (F_1)$}
 \put(200,60){\circle*{3}}
 \put(210,70){\circle*{3}}
 \put(150,70){\circle*{3}}
 \put(150,80){\circle*{3}}
 \put(152,80){\vector(1,0){6}}
 \put(150,82){\vector(0,1){6}}
 \put(150,90){\circle*{3}}
 \put(160,80){\circle*{3}}
 \put(198,156){\line(1,0){4}}
\put(200,60){\line(1,1){90}} \put(295,155){\line(1,1){10}}
\put(200,60){\line(-5,1){80}} \put(120,76){\line(1,1){110}}
\put(235,191){\line(1,1){10}} \put(170,90){
\rotatebox{45}
{$K^{\bullet,\bullet}_i$ nonzero only in this
strip}}
\end{picture}
\caption{\label{fig:K}the double complex $K^{\bullet,\bullet}_i$.}
\end{figure}

%Figure 1 illustrates that $K_i^{\bullet,\bullet}$
%is nonzero only if
%$0\le -p+q\le \rk(F_1)$ and $p+mq\ge 0$.

Note that the direct sum of the
double complexes $K_i^{\bullet,\bullet}$ is equal to
$L^{\bullet,\bullet}$.
$$\bigoplus_{i=0}^m K_i^{\bullet,\bullet}=L^{\bullet,\bullet}$$
For any $i=0,\dots,m$, we denote by $(K_i^\bullet, D=\wt d+\wt a)$
the total complex.
%The complex $K_i^{\bullet,\bullet}$ is nonzero only if
%$0\le -p+q\le \rk(F_1)$ and $p+mq\ge 0$.
\begin{thm}\label{bound}
There exists an integer $n_0$ such that,
for  $n\ge n_0$,
the cohomology group ${H}^n(K_i^\bullet)$
vanishes and
the sheaf
$I=\im(K_i^{n-1}\to K_i^{n})$ is
locally free.
\end{thm}

\begin{proof}
We write $K^{\bullet,\bullet}$
and $K^\bullet$ omitting $i$.
We define
\begin{equation}\label{n_0}
n_0=(4m+2\rk(\Sym^m F_0))/(m+1)
\end{equation}
and show the exactness of $K^{n-1}\to K^n\to K^{n+1}$
for $n\ge n_0$.

Since all the sheaves of $K^\bullet$
are coherent and locally free,
it suffices to check the exactness after base change
of $K^\bullet$ to $x$ for all closed points
$x\in X$.
The spectral sequence
of complex vector spaces
${E}_2^{p,q}\cong { H}^q_{\widetilde a} ( {
H}^p_{\widetilde d}(K^{\bullet,\bullet}_x))$ abuts to the cohomology
of $(K^\bullet_x,D_x)$. We show $E_2^{p,q}=0$ for
$p+q\ge n_0$.

Write $F_{1,x}\to F_{0,x}$ as
the sum of
$H_{1,x}\xrightarrow{0} H_{0,x}$ and
$V\xrightarrow{\id} V$.
In this way, each term $K^{p,q}_x$ can be written
as the sum of the complex vector spaces $\Sym^{r} H_{0,x}\otimes
\Lambda^{s}H_{1,x}\otimes \Sym ^{p+mq-i-r} V\otimes
\Lambda^{-p+q-s} V$ for $r\in\{0,\dots, p+mq-i\}$ and $s\in \{0,\dots, -p+q\}$.
Denote by $G^{p,q}$ the summand corresponding to $r=0$ and $s=0$
and denote by
$Q^{p,q}$ the sum of the remaining terms in $K^{p,q}_x$.
In this way, each element
$\al \in K^{p,q}_x$ is given by $\al'\in G^{p,q}$
and $\al''\in Q^{p,q}$.
The morphism
$a_x\colon \CC_x\to \Sym^m F_{0,x}\otimes F_{1,x}= K^{0,1}$
splits as $(a'_x\in G^{0,1},a''_x\in Q^{0,1})$.
The differentials $\wt d_x$ and
$\wt a_x$  can be written
as
\begin{align}
G^{p,q}\oplus Q^{p,q}\to G^{p+1,q}\oplus Q^{p+1,q} && \wt d_x&=
\begin{pmatrix}0&0\\0 & d''\end{pmatrix}, \\
G^{p,q}\oplus Q^{p,q}\to G^{p,q+1}\oplus Q^{p,q+1} && \wt a_x&=
\begin{pmatrix}\ol a'&0\\\ol a'' & \widehat{a}\end{pmatrix},
\label{splita}
\end{align}
where, by construction, $d''$ is exact, and
we have $d''\circ\widehat{a}+\widehat{a}\circ d''=0$
and $d''\circ \ol a''=0$, as a consequence of
$\wt d_x\circ \wt a_x+\wt a_x\circ \wt d_x=0$,
(Lemma \ref{doublecx}).

Now, note that
$H^p_{\widetilde d}(K_x^{\bullet,q})=G^{p,q}$,
by the exactness of $d''$.
Furthermore, the vertical cohomology
$H^q_{\widetilde a}(H^p_{\widetilde d}(K_x^{\bullet,\bullet}))$
is the cohomology of the complex
$G^{p,\bullet}$ where the differential $a'$ is
induced by $a'_x\colon \CC_x\to G^{0,1}=
\Sym^m H_{0,x}\otimes H_{1,x}$ as follows
\begin{multline}\label{ola}
\Sym^h H_{0,x}\otimes \Lambda^k H_{1,x}
\xrightarrow{\ \id\otimes a'_x \otimes \id \  }
\Sym^h H_{0,x}\otimes G^{0,1}
\otimes \Lambda^k H_{1,x}\\
=\Sym^h H_{0,x}\otimes \Sym^m H_{0,x}\otimes H_{1,x}
\otimes \Lambda^k H_{1,x}\\
\xrightarrow{\ \sigma_m\otimes \lambda_1 \ }
\Sym^{h+m} H_{0,x}\otimes \Lambda^{k+1} H_{1,x}.
\end{multline}
Indeed, by \eqref{splita}, 
we have $\wt a_x(\al',\al'')=
(\ol a'(\al'),\ol a''(\al')+\widehat{a}(\al''))$.
By the exactness of $d''$, the passage to 
cohomology with respect to
$\wt d_x$ induces
the homomorphism
$\ol a'$ of \eqref{ola} on the cohomology
groups $H^p_{\widetilde d}(K_x^{\bullet,q})$.
So, $$H^p_{\widetilde d}(K_x^{\bullet,q-1})\to
H^p_{\widetilde d}(K_x^{\bullet,q})
\to H^p_{\widetilde d}(K_x^{\bullet,q+1}),$$
is the
cohomological Koszul complex
associated to $(H_{1,x})^\vee\to \Sym^{m}H_{0,x}$:
\begin{multline}
\label{vertcx}\Sym^{h-m}H_{0,x}\otimes\Lambda^{k-1}H_{1,x}
\to \Sym^{h}H_{0,x}\otimes\Lambda^{k}H_{1,x}\\
\to \Sym^{h+m}H_{0,x}\otimes\Lambda^{k+1}H_{1,x},
\end{multline}
for $h=p+mq-i$ and $k=-p+q$.
After tensoring by $\det (H_{1,x})^\vee$, we
obtain the homological Koszul differential
\begin{multline}
\label{t} {\Sym}^{w-m}H_{0,x}\otimes\Lambda^{v+1}H_{1,x}^\vee
\rightarrow {\Sym}^{w}H_{0,x}\otimes
\Lambda^{v}H_{1,x}^\vee\\
\rightarrow
{\Sym}^{w+m}H_{0,x}\otimes \Lambda^{v-1}H_{1,x}^\vee,
\end{multline}
with $v=\rk(H_{1,x})-q+p$ and $w=p+mq-i$.

We illustrate
that the homomorphism
of \eqref{t} can be regarded as 
the homological Koszul complex associated 
to the  linear system
\begin{equation}\label{map}
S_x\colon H_{1,x}^\vee\to
{\Sym}^m H_{0,x}=H^0(Y,{\mathcal
O}_Y(m)),
\end{equation}
where  $Y$ is $\PP H_{0,x}^\vee$.
Indeed, denote by $J_x$
the image of $H_{1,x}^\vee$ via $S_x$; recall that
$J_x$
is a base point free linear system
because $a$ is nondegenerate.
Write \linebreak $H_{1,x}^\vee=N_x \oplus J_x$,
where $N_x$ is the kernel
of $S_x$. The composite homomorphism at (\ref{t})
is the sum over $t\in \{0,\dots, v+1\}$ of
\begin{multline}\label{N&J}
{\Sym}^{w-m}H_{0,x}\otimes
\Lambda^{v-t+1}J_x
%\otimes \Lambda^t N_x
\rightarrow
{\Sym}^{w}H_{0,x}\otimes
\Lambda^{v-t}J_x
%\Lambda^t N_x
\\
\rightarrow
{\Sym}^{w+m}H_{0,x}\otimes
\Lambda^{v-t-1}J_x
%\otimes \Lambda^t N_x
\end{multline}
tensored by $\Lambda^t N_x$.

In \cite{Gr}, Green treats the
case of the homological
Koszul complexes
induced by a base point free linear system.
By \cite[Thm.~2]{Gr}, the sequence
\eqref{N&J} is exact
when $J_x$ is base point free
and the condition
$$w\ge (v-t)+m+\codim(J_x)=v-t+m+\rk({\Sym}^m H_{0,x})-\rk(J_x)$$
is satisfied.
The middle term of \eqref{N&J}
is nonzero only if
$v-t-\rk(J_x)\le 0$.
Therefore, \eqref{N&J} is exact for
any $t$ if
$w\ge m+ \rk ({\Sym}^m H_{0,x})$.

Recall that $K^{p,q}\ne 0$ only if
$q-p\ge 0$; therefore,
we restrict to the indexes $(p,q)$
satisfying $q-p\ge 0$.
In this way, we have
$(m-1)(q-p)\ge 0$. Furthermore, we have 
$(m+1)(p+q)\ge (m+1)n_0$, which follows from
$p+q\ge n_0$.
Summing up, we get
$2p+2mq\ge (m+1)n_0$ and, by \eqref{n_0},
$2p+2mq\ge 4m +2 \rk (\Sym ^m F_0)$, which
implies
$$p+mq-m\ge m + \rk (\Sym ^m F_0).$$
We obtain
$$w=p+mq -i\ge p+mq-m\ge m+ \rk (\Sym ^m F_0)
\ge m + \rk (\Sym ^m H_{0,x});$$
therefore, \eqref{N&J} is exact and $E_2^{p,q}$ vanishes.

Finally,
since $H^{n}((K_i^\bullet)_x)$ vanishes for
every $x\in X$ and for $n\ge n_0$,
we have, for all $n\ge n_0$,
$$\ker(D_{n})_x = \im(D_{n-1})_x,$$
where $(D_n)_x$ is the homomorphism
induced by the total differential \linebreak $D_n\colon K^n\to K^{n+1}$
on the fibre over $x$.
In fact, the rank
of $\ker(D_{n})_x$ is equal to the rank of $\im(D_{n-1})_x$
and is constant, because it is upper and lower semicontinuous.
Therefore, the sheaf
$I=\im(K_i^{n-1}\to K_i^{n})$ is
locally free.
\end{proof}
\subsection{The $K$-class.}
By Theorem \ref{bound},
the alternate sum of the
cohomology sheaves of $K_i^\bullet$
is a well defined class in the $K$-theory
group $K_0'(X)$.
The fact
that $I=\im(K_i^{n-1}\to K_i^{n})$ is
locally free for $n\ge n_0$ allows us to define the
$K$-class
$$\class(K^\bullet_i)=
\sum_{n<n_0} (-1)^n[K_i^n]+(-1)^{n_0}[{ I}]\in K_0(X).$$
which lifts the natural class $\sum_n(-1)^n[{
H}^n(K^\bullet_i)]$ in $K_0'(X)$.

\begin{defn}\label{keulerclass}
Let $F_\bullet$ be a complex in degrees $1$ and $0$
with a closed and
nondegenerate form $a$. Then, the
{\em $K$-theory Euler class} of
$(F_\bullet,a)$
is $$\ke(F_\bullet,a)= \sum_{0\le
i\le m}\class(K_i^\bullet).$$
\end{defn}
\begin{pro}\label{pro:oneiszero}
Let $V$ be a vector bundle of finite rank on $X$.
Consider the complex $V\to 0$
in degrees $1$ and $0$.
A closed and nondegenerate
form $a$ is necessarily the zero
homomorphism $\mathcal O\to 0$.
The $K$-theory Euler class extends
the total lambda class evaluated at $-1$
in the sense that
\begin{equation}\label{kevect}
\ke(V\to 0, a=0)=\lambda_{-1}(V).
\end{equation}
\end{pro}
\begin{proof}
By Definition \ref{keulerclass}, we have
\begin{multline*}\ke(V\to 0, 0)=\sum_{0\le i\le m}
\class(K^\bullet_i)=\sum_{0\le i\le m}
\sum_{\substack{
p, q\ge 0 \\
p+mq=i}}
(-1)^{q-p}[\Lambda^{q-p}V]\\
=\sum_{0\le i\le m}
\sum_{\substack{
q\ge 1}}(-1)^{(m+1)q-i}[\Lambda^{(m+1)q-i}V]
=
\sum_{\substack{q\ge 0}}(-1)^{q}[\Lambda^{q}V].
\end{multline*}
\end{proof}

\subsection{The cohomology class.}
We define the cohomological realization
of $\ke(F_\bullet,a)$ by
analogy with \eqref{ke-ctop}.
\begin{defn}\label{defnctop}
Let $E^0\to E^1$ be a complex of coherent
and locally free sheaves
in degrees $0$ and $1$,
and $b\colon \Sym^mE^0\otimes E^1\to \mathcal O_X$ a
homomorphism such that
$b^\vee$ is a closed and nondegenerate
form for $(E^\bullet)^\vee$.
Then, $$\ctop(E^\bullet,b)=
{\ch}({\ke}((E^\bullet)^\vee,b^\vee))\frac{{\td}(E^1)}{\td(E^0)}$$
is the {\em top Chern class} of $(E^\bullet,b)$
in the rational Chow ring $A^*(X)_\QQ$.
\end{defn}
\begin{rem}\label{comp-w-ctop}
Clearly, by Proposition \ref{pro:oneiszero} and \eqref{ke-ctop}, we have
$$\ctop(0\to V, b=0)=\ctop(V).$$
\end{rem}

\subsection{Dependence on $F_\bullet$ and $a$.}
We investigate how $\ke$
depends on the complex $F_\bullet$ and on the
form $a$.
\begin{thm}\label{indep}
The class $\ke(F_\bullet,a)$ depends on $F_\bullet$ up to
quasiisomorphism and on the form $a$ up to homotopy.
\end{thm}
\begin{proof}
We prove the claim in two lemmata.
\begin{lem}\label{homa}
Let $a$ and $a'$ be two closed and nondegenerate forms
\linebreak $\mathcal O_X\to \Sym^m F_0\otimes F_1$ defining
homotopic homomorphisms of complexes:
\[ \xymatrix{
\dots\ar[r]&
    \Sym^{m-1}F_0\otimes \Lambda^2 F_1\ar[r]&
            \Sym ^{m}F_0\otimes F_1 \ar[r]&
                \Sym^{m+1} F_0
                                \\
   \dots \ar[r]&
        0\ar[r]\ar[u] &
            {\mathcal O}_X \ar[u]_{a,a'}\ar[r]
                & 0.\ar[u]              } \]
Then, for any $i=0,\dots,m$, the total complex
$(K_i^\bullet,D={\widetilde d}+{\widetilde a})$
and the total complex $(K_i^\bullet,D'={\widetilde d}+{\widetilde a}')$
are isomorphic. In particular, in $K$-theory, we have the identity
$\ke(F_\bullet,a)=\ke(F_\bullet,a')$.
\end{lem}
\begin{proof}
Let $h$ be the
homotopy ${\mathcal O}_X\to\Sym^{m-1}F_0\otimes \Lambda^2 F_1$
making the \linebreak diagram
\[ \xymatrix@C=1cm{
\Sym^{m-1}F_0\otimes \Lambda^2F_1\ar[r]^{ d}&
    \Sym^{m}F_0\otimes F_1 \\
   &
    {\mathcal O}_X\ar[ul]^{h}\ar[u]_{a - {a}'} \\
} \] commutative:
\begin{equation}\label{h}
{d}\circ {h}={a}-{a}'.
\end{equation}
In fact, $h$ allows us to define a map $\widetilde h\colon
K_i^{p,q}\rightarrow K_i^{p-1,q+1}$ as the natural composite
homomorphism
\begin{multline*}{\Sym}^{p+mq-i}F_0\otimes{\mathcal O}_X\otimes\Lambda^{q-p}F_1
\\ \rightarrow{\Sym}^{p+mq-i}F_0 \otimes
{\Sym}^{m-1}F_0\otimes\Lambda^2F_1 \otimes\Lambda^{q-p}F_1
\\\rightarrow
{\Sym}^{p-1+m(q+1)-i}F_0\otimes\Lambda^{q-p+2}F_1.\end{multline*}
More explicitly, if we write $h$
locally as
$\sum_l H_l \otimes \sigma_{l}^+\wedge  \sigma_{l}^-$,
the homomorphism $\widetilde h$ sends
$\prod_{1\le i\le h}c_i\otimes \bigwedge_{1\le j\le k} b_j$
to \begin{equation}
\label{honfibre} \sum_{l}\prod_{1\le i\le h}c_i\cdot
H_l\otimes \sigma_{l}^+\wedge\sigma_{l}^-\wedge
\bigwedge_{1\le j\le k} b_j.
\end{equation}

Now, by means of $\widetilde h$ we define an isomorphism
between the two double complexes.
First, we need to show two relations:
\eqref{tildeh} and \eqref{dhcommk}.

By \eqref{h}, we have immediately
\begin{equation}\label{tildeh}
\widetilde a -\widetilde a'=\widetilde d \circ \widetilde h -
\widetilde h \circ \widetilde d.
\end{equation}
Furthermore, we have
\begin{equation}\label{dhcommute}
\widetilde d\circ \widetilde h\circ \widetilde h
=- \widetilde h\circ \widetilde h\circ \widetilde d
+ 2 \widetilde h\circ \widetilde d\circ \widetilde h
\end{equation}
as a straightforward consequence of
$\widetilde d( \widetilde h(h))=
2\widetilde h(\widetilde d(h))$, which
we can verify by means of
the local presentation of $\widetilde h$
given above. Indeed, apply $\widetilde h$ to
$$\widetilde d(h)= \sum_l H_l\cdot d(\sigma_{l}^+)\otimes
\sigma_{l}^- - \sum_l H_l\cdot \sigma_{l}^-\otimes
d(\sigma_{l}^+)$$
and compare with the image via $\widetilde d$ of
$$\widetilde h(h)=
2\sum_{l,q} H_l\cdot H_q\otimes \sigma_{l}^+\wedge
\sigma_{l}^-\wedge \sigma_{q}^+\wedge
\sigma_{q}^-.$$

Using \eqref{dhcommute}, we prove by induction
the following equation for $n\ge 0$
\begin{equation}\label{dhcommk}
n(\widetilde d\circ \widetilde h^{n+1})
=- \widetilde h^{n+1} \circ \widetilde d
+ (n+1)( \widetilde h\circ
\widetilde d\circ \widetilde h^{n}),
\end{equation}
where $\widetilde h^n$ is the
composition of $\widetilde h$
iterated $n$ times.
First, applying the equation for $n-1$,
we see that the right hand side equals
\begin{multline*}
(- \widetilde h^{n+1} \circ \widetilde d
+ \widetilde h\circ \widetilde d\circ \widetilde h^{n})
+ n( \widetilde h\circ \widetilde d\circ \widetilde h^{n-1})
\circ \widetilde h\\
=(- \widetilde h^{n+1} \circ \widetilde d
+ \widetilde h\circ \widetilde d\circ \widetilde h^{n})
+ (\widetilde h^{n} \circ \widetilde d
+ (n-1)(\widetilde d\circ \widetilde h^{n}))\circ \widetilde h.
\end{multline*}
By \eqref{dhcommute}, for all
integers $0\le i\le n+1$, we have
$-\widetilde h^i\circ \widetilde d\circ \widetilde h^{n+1-i}
+\widetilde h^{i+1}\circ \widetilde d\circ \widetilde h^{n-i}=
-\widetilde h^{i+1}\circ \widetilde d\circ \widetilde h^{n-i}
+\widetilde h^{i+2}\circ \widetilde d\circ \widetilde h^{n-i-1}$,
which implies \eqref{dhcommk} by iterated application.

Now we introduce a homomorphism of total complexes:
$$e^{-\widetilde h}:=\sum_{j\ge 0} (-\widetilde h)^{j}/j! \ \ .$$
Since $\widetilde h$ respects the total
grading, it is an endomorphism of $K^n_i$.
It is well defined (in characteristic 0)
because the summands above vanish
as soon as $j> (\rk F_1)/2$
(this immediately follows from \eqref{honfibre}).
For all $n$, $e^{-\widetilde h}$ is
in fact an automorphism of $K_i^n$, because
$e^{\widetilde h}$ is also well defined
and is the
inverse of $e^{-\widetilde h}$.
Finally,
putting together \eqref{tildeh} and \eqref{dhcommk},
we prove
$$e^{-\widetilde h} \circ D=D'\circ e^{-\widetilde h}.$$
Indeed, by \eqref{tildeh}, we have
$D'=\widetilde d+
\widetilde a +\widetilde d\circ \widetilde h
- \widetilde h\circ \widetilde d$, and it is enough to
show
$$e^{-\widetilde h}\circ \widetilde d
=\widetilde d\circ e^{-\widetilde h}+
 \widetilde d\circ \widetilde h\circ e^{-\widetilde h}
+\widetilde h \circ \widetilde d\circ e^{-\widetilde h},$$
because $\widetilde h$ and $\widetilde a$
clearly commute.
Indeed, in the equation above,
the terms of  degree $n$ in $\widetilde h$
yield precisely the equation \eqref{dhcommk}.
\end{proof}

\begin{lem}\label{qi}
Let $(F_\bullet,a)$ and $(\Phi_\bullet,\al)$ be two pairs
formed by a complex and a closed and nondegenerate form.
Let $f\colon \Phi_\bullet\to F_\bullet$ be
a quasiisomorphism with
\begin{equation}\label{atoal}
a=\al\circ \Lambda^{m+1}(f).
\end{equation}
Let $\Gamma_i^{\bullet,\bullet}$ and $K_i^{\bullet,\bullet}$ be the
double complexes associated to $(\Phi_\bullet,\al)$ and
$(F_\bullet,a)$. Then, the total complexes
$K_i^\bullet$ and
$\Gamma_i^\bullet$ are
quasiisomorphic. In particular,
we have $\ke(F_\bullet,a)=\ke(\Phi_\bullet,\al).$
\end{lem}
\begin{proof}
We denote by $\wt a$ and $\wt d$ the differentials of
$K_i^\bullet$ as in Section \ref{doubleK}.
We denote by $\widetilde \delta$ and $\widetilde \al$
the analog differentials of bidegrees $(1,0)$ and $(0,1)$
of $\Gamma_i^{\bullet,\bullet}$.
Define $\wt f\colon \Gamma_i^{\bullet,\bullet}\to
K_i^{\bullet,\bullet}$ by the equation
$$\wt f=\Lambda^{mq+q-i}(f):\Sym^{p+mq-i}\Phi_0
\otimes \Lambda^{-p+q}\Phi_1
\rightarrow
\Sym^{p+mq-i}{F}_0\otimes \Lambda^{-p+q}{F}_1.$$
Note that
$\wt f$ satisfies
$\wt f\circ \wt \delta= \wt d \circ \wt f$ and, by
(\ref{atoal}),  $\wt f\circ \wt \alpha= \wt a \circ \wt f$;
therefore, $\wt f$ is a homomorphism of
double complexes which respects the filtrations
$\sum_{q\ge t}K_i^{p,q}$
and $\sum_{q\ge t}{\Gamma}_i^{p,q}$.

Note that $\wt f$ is a morphism
of filtered differential
graded modules and
induces homomorphisms of bigraded
modules $\phi_l: E_l^{p,q}(K^{\bullet,\bullet})\to
E_l^{p,q}(\Gamma^{\bullet,\bullet})$ for $l\ge 0$.
Note that $\phi_1$ is
an isomorphism, because
$\Lambda^{N}(f)\colon
\Lambda^N(\Phi_\bullet)\to \Lambda^N(F_\bullet)$
is a quasiisomorphism.
Then, by \cite[\S 3.1, Thm.~3.2]{Mc01},
$\phi_l$ is an isomorphism
for every $l\ge 1$.
\end{proof}
The theorem follows from Lemma \ref{homa} and Lemma \ref{qi}.
\end{proof}
\subsection{Passage to the derived category.}
\label{sect:derived}
In the next section we consider
a complex $E^\bullet$
in degrees $0$ and $1$ of coherent locally free sheaves on $X$
and a morphism
in the derived category
${\bf D}^b(X)$
$$\tau \colon \Sym^{m+1}(E^\bullet)\to
\mathcal O_X[-1]$$
inducing a base point free linear system
$S_x\colon H^1(E^\bullet_x)\rightarrow
\Sym^m H^0(E^{\bullet}_x)^\vee$ for any $x\in X$.
When $X$ is quasiprojective,
Definition \ref{keulerclass}
and Definition \ref{defnctop}
descend to the derived category
in the same way as in \cite{PV}.

Indeed, by \cite[Prop.~4.7]{PV}, for $h$ and $k$ sufficiently large,
there exists a complex
$C^\bullet=(C^0\to \mathcal O_X(-h)^{\oplus k})$
of coherent locally free sheaves on $X$,
quasiisomorphic to $E^\bullet$, and satisfying
\begin{equation}\label{liftPV}
\Hom_{{\bf K}^b(X)}(\Sym^{m+1}(C^\bullet),
\mathcal O_X[-1])\cong
\Hom_{{\bf D}^b(X)}(\Sym^{m+1}(C^\bullet),
\mathcal O_X[-1]),
\end{equation}
where ${\bf K}^b(X)$ is the
homotopic category of bounded
complexes of coherent sheaves on $X$.
Therefore,
the  morphism $\tau$
in the derived category
is lifted up to homotopy
by the homomorphism of complexes
$\gamma_\bullet\colon
\Sym^{m+1}(C^\bullet)\to
\mathcal O_X[-1]$.

Let $c$ be the homomorphism $\gamma_1\colon
\Sym^m C^0\otimes C^1 \to\mathcal O_X$.
Note that $c^\vee$ is a closed and nondegenerate form
for $(C^\bullet)^\vee$; therefore,
we define
\begin{align}
\ke((E^\bullet)^\vee,\tau^\vee)=\ke((C^\bullet)^\vee,c^\vee).
\end{align}
Furthermore, Definition \ref{defnctop} applies to
$C^\bullet$ and $c$. We define
\begin{equation}
\ctop(E^\bullet,\tau)=\ctop(C^\bullet,c).
\end{equation}
The definitions above do not depend on the
choice of $C^\bullet$ and $c$.
Indeed, by
\cite[Lem.~4.8]{PV},
for any quasiisomorphisms
$C^\bullet\to E^\bullet$ and $D^\bullet\to E^\bullet$
of complexes in degrees $0$ and $1$
and for any homomorphisms lifting $\tau$
\begin{align*}
c\colon &\Sym^m C^0\otimes C^1 \to \mathcal O_X, &
d\colon & \Sym^m D^0\otimes D^1 \to \mathcal O_X,
\end{align*}
there exists a complex
$\ol E^\bullet=(\ol E^0\to
\mathcal O_X(-h)^{\oplus k})$,
and two quasiisomorphisms
$f\colon \ol E^\bullet\to C^\bullet$,
$g\colon \ol E^\bullet\to D^\bullet$.
By \cite[Prop.~4.7]{PV}, for a sufficiently large $h$,
we can choose $\ol E^\bullet$ in such a way
that the isomorphism \eqref{liftPV} holds.
The homomorphisms
$\ol c=c \circ \Sym^{m+1}(f)$
and $\ol d=d\circ \Sym^{m+1}(g)$
lift $\tau$. Therefore, by \eqref{liftPV},
they are homotopic.
Then, by Lemma \ref{qi}, we have
$\ke((C^\bullet)^\vee,c^\vee)=\ke((\ol E^\bullet)^\vee,\ol c^\vee)$ and
$\ke((D^\bullet)^\vee,d^\vee)=\ke((\ol E^\bullet)^\vee,\ol d^\vee)$,
whereas, by Lemma \ref{homa}, we have
$\ke((\ol E^\bullet)^\vee,\ol c^\vee)=\ke((\ol E^\bullet)^\vee,\ol d^\vee)$.

\begin{rem}
The $K$-class $\ke(F^\bullet,a)$
of Definition \ref{keulerclass}
only depends on $F_\bullet$ and
on the homomorphism $\tau$ induced by $a$
in the derived category.
The $K$-class can only be defined
when there exists a suitable
homomorphism $\tau$. It would be interesting
to determine more accurately
how it depends on it.
\end{rem}

\section{The Witten top Chern class}\label{sect:wtcc}
We choose nonnegative integers $r$, $g$, $n$, and
$(k_1,\dots,k_n)={\kb}$
satisfying the conditions
$r\ge 2$, $2g-2+n>0$, and $2g-2-\sum k_i\in r\ZZ$.
Let $\MSbar$ be the moduli stack of
stable $r$-spin curves.
By \cite{Jageom}, an object of
$\MSbar$ is a set of data including:
\begin{enumerate}
\item a stable curve
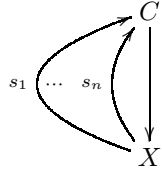
\[ \xymatrix@R=1.5cm{
{C}\ar[d]\\
X\ar@/^1.5cm/[u]^{s_1}_\cdots\ar@/^0.5cm/[u]^{s_n}\\
} \]
where $\pi$ is proper,
the sections $s_i$ have disjoint images in
the smooth locus of $\pi$,
and $X$ is a
quasiprojective scheme,
\item
a torsion free sheaf of rank one $L$ on $C$,
\item
an injective homomorphism
$$f\colon L^{\otimes r}\hookrightarrow
\omega_{C/X} \left(-\sum\nolimits_i k_i[
s_i(X)]\right).$$
\end{enumerate}

The data
$(\pi\colon C\to X;L;f)$ determine
$\tau \colon \Sym^r(R\pi_*L)\to
{\mathcal O}_X[-1]$ in the derived category ${\bf D}^b(X)$.
Indeed, by \cite[Sect.~5]{PV},
the homomorphism $\Sym^r(R\pi_*L)\to
{\mathcal
O}_X[-1]$ is the composite of
the trace homomorphism
from Grothendieck
duality  $R\pi_*\omega_{{C}/X}\to
{\mathcal O}_X[-1]$ and the homomorphism
$\Sym^r(E^\bullet)\to R\pi_*\omega_{C/X}$ induced by
\begin{equation}\label{gamma}
\gamma\colon L^{\otimes r}\xrightarrow{f}
\omega_{C/X} \left(-\sum\nolimits_i k_i[
s_i(X)]\right)\to\omega_{C/X}.\end{equation}

It is straightforward that
the class  $R\pi_*L$ can
be represented by a complex $E^\bullet$ in degrees $0$ and $1$.
On the other hand, for every $x\in X$,
the homomorphism $\tau$ induces
\begin{equation}\label{pairing}
\Sym^{r-1}H^0(E^\bullet_x) \otimes H^1(E^\bullet_x)\to \CC,
\end{equation}
where $H^i(E^\bullet_x)$ is isomorphic to
$H^i(C_x,L_x)$.
The homomorphism \eqref{pairing} \linebreak 
coincides with the composite homomorphism
\begin{multline}\label{tau}
\tau_x\colon\Sym^{r-1}H^0(E^\bullet_x)\rightarrow
H^0\left({C}_x,L_x^{\otimes r-1}\right)\hookrightarrow
H^0\left({C}_x,\omega_{{C}_x}\otimes {L_x}^\vee\right)\\
\xrightarrow{\sim}
(H^1(E^\bullet_x))^\vee.
\end{multline}
Note that an element $v\in H^0(E^\bullet_x)$
is zero if for every $w\in H^1(E^\bullet_x)$ we have
$\langle\tau_x(v^{r-1}),w\rangle=0$
(the injectivity of the middle homomorphism in \eqref{tau}
is crucial and is guaranteed by $k_i\ge0$).
This means that $\tau$ induces a base point
free linear system; therefore,
choosing $m=r-1$,
we get a well defined $K$-class
$\ke((E^\bullet)^\vee,\tau^\vee)\in K_0(X)$.

For any object of $\MSbar(X)$,
we define the classes
$\ke((E^\bullet)^\vee,\tau^\vee)$ and \linebreak
$\ctop(E^\bullet,\tau)$ in
$K_0(X)$ and $A^*(X)_\QQ$.
The constructions are compatible with
morphism in the category
$\MSbar$ and naturally induce a $K$-theory class
in $K_0(\MSbar)$ and a Chow cohomology class in
$A^*(\MSbar)_\QQ$.
\begin{defn}\label{wtcc-chow}
We define the
\emph{Witten top Chern class in $K$-theory} as
$$K_W=\ke((E^\bullet)^\vee,\tau^\vee)$$ and
the \emph{Witten top Chern class} as
\begin{equation}
c_W=\ctop(E^\bullet,\tau).
\end{equation}
\end{defn}

\begin{rem}[Witten top Chern class using other compactifications]
\label{compstack}
In this paper we work with Jarvis's compactification \cite{Jageom}.
In fact, 
one can give a different
compactification of the stack of smooth $r$-spin curves $\MS$
by means of (balanced) twisted curves
in the sense of Abramovich and Vistoli.
In this way, for instance, Abramovich and Jarvis rephrase the
original construction of Jarvis,
\cite{AJ}.
A different approach exploits
Olsson's description of the category of twisted curves
\cite{Ol}:
in \cite{Chmod},
we consider the stack $\widetilde {\stack S}_g(r)$
of $r$-spin structures
over twisted curves,
which is \'etale over Olsson's
stack of twisted curves $\widetilde \MM_g$:
$$\widetilde {\stack S}_g(r)\to \widetilde {\MM}_g$$
(in sketching the procedure, we are omitting
the markings for simplicity).
We proceed by classifying all the
compactifications of $\MM_g$ contained
in $\widetilde \MM_g$. Then,
base change leads to new compactifications.

Although the new compactifications are not isomorphic
to the preexisting ones, the rational Chow rings are
isomorphic and
the construction can be applied without modifications
and yields the same class in rational cohomology.
This happens because, the pushforward
homomorphism $R(\pi_{\stack C})_*$ of
the $r$-spin structure $\stack L$
on the twisted curve $\pi_{\stack C}\colon \stack C\to X$
in the derived category
factors through the derived category of the
coarse space $C$. Now, for stacks of
Deligne--Mumford type,
the pushforward homomorphism via $\stack C\to C$
is exact on coherent sheaves \cite[Lem.~2.3.4]{AV}.
Finally, this procedure yields the class $c_W$,
because it carries $r$-spin structures
on twisted curves
to Jarvis's $r$-spin structures on stable curves defined using
relatively torsion-free
sheaves (see \cite[\S3, \S4.3]{AJ}).
\end{rem}

\section{Compatibility with the Polishchuk--Vaintrob class}\label{sect:pv}
\subsection{Notation.}
The key ingredient of the construction is a $2$-periodic
complex.
\begin{defn} A {\em $2$-periodic complex} of sheaves on a scheme
$Y$ is a
$\ZZ/2\ZZ$-graded sheaf $W=W^+\oplus W^-$ on
$Y$ equipped with $d^+\colon W^+\to W^-$ and $d^-\colon W^-\to
W^+$ with $d^-\circ d^+=d^+\circ d^-=0$. We say that $W$ is a
complex of coherent (respectively quasicoherent) sheaves if $W^+$
and $W^-$ are coherent (respectively quasicoherent). We say that
$W$ is exact if
$$\cdots\xrightarrow{d^-} W^+\xrightarrow{d^+}W^-\xrightarrow{d^-}
W^+\xrightarrow{d^+}W^-\xrightarrow{d^-}\cdots$$ is exact.
\end{defn}

We denote by $\EE=\Spec (\Sym^*E^\vee)$ the total space of a
coherent locally free sheaf $E$ on a scheme $X$.

\subsection{The construction of Polishchuk and Vaintrob.}
\label{par:pv}
We recall the construction from \cite{PV}.
We choose an integer $m\ge 1$.
Consider  a complex $E^\bullet$
in degrees $0$ and $1$
of coherent locally free sheaves on $X$
and a homomorphism $b\colon \Sym^m E^0\otimes E^1 \to
\mathcal O_X$ such that $b^\vee$ is
a closed and nondegenerate form for $(E^\bullet)^\vee$.

Let $\EE^0$ be the total space of $E^0$,
$p$ the
projection $\EE^0\to X$, and
$i\colon X\to \EE^0$ the zero section.
The pullback of  $E^0\to E^1$ via $p$ corresponds to a section
$\delta\colon {\mathcal O}_{\EE^0}\to p^*E^1$.
Compose the natural homomorphism $E^0\to
\Sym^{m}E^0$ with the homomorphism
${\Sym}^{m} E^0\rightarrow ({E^1})^\vee$
induced by $b$. The pullback via $p$
corresponds to a section $\al\colon {\mathcal
O}_{\EE^0}\to p^*(E^1)^\vee$.

Denote by $S^h$ the sheaf $\Lambda^{h}p^*(E^1)^\vee$ on $\EE^0$;
then, $\delta$ and $\al$ allow us to define homomorphisms
$\widetilde \delta\colon S^{h+1}\to S^h$ and $ \widetilde
\al\colon S^{h}\to S^{h+1}$. The constructions of $\wt \delta$
and of $\wt \al$ are analogous to the ones of $\wt d$ and $\wt
a$. For example, to define $\widetilde \delta$, we take the
natural map $\Lambda^{h+1}p^*(E^1)^\vee\to
\Lambda^{h}p^*(E^1)^\vee \otimes p^*(E^1)^\vee$, the obvious
pairing $p^*(E^1)^\vee\otimes p^*(E^1)\to {\mathcal O}_{\EE^0}$,
and the composite homomorphism
\begin{multline*}
\widetilde \delta\colon \Lambda^{h+1}p^*(E^1)^\vee\to
\Lambda^{h}p^*(E^1)^\vee \otimes p^*(E^1)^\vee
\otimes {\mathcal O}_{\EE^0}\\
\xrightarrow{\id\otimes\id\otimes\delta} \Lambda^{h}p^*(E^1)^\vee
\otimes p^*(E^1)^\vee \otimes  p^*E^1\to \Lambda^{h}p^*(E^1)^\vee.
\end{multline*}
Define $d^+\colon S^+\to S^-$ and $d^-\colon S^-\to S^+$ by
summing for $h$ even and $h$ odd the homomorphisms
$\widetilde\delta\colon S^h\to S^{h-1}$ and $\widetilde \al\colon
S^{h}\to S^{h+1}$.
It is shown in \cite{PV} that the homomorphisms
$d^+\circ d^-$ and $d^-\circ d^+$ vanish. We denote by
$S=S^+\oplus S^-$ the $2$-periodic complex of coherent locally
free sheaves on $\EE^0$ with differential $d=d^+\oplus d^-$.
(In fact $S$ is the spinor bundle associated to the
orthogonal bundle $p^*E^1\oplus p^*(E^1)^\vee$ and the
differential is naturally induced by the isotropic section
$\delta\oplus\al$.)

The complex
\begin{equation}\label{P-cx}
\cdots\xrightarrow{d^-}
S^+\xrightarrow{d^+}
S^-\xrightarrow{d^-} S^+\xrightarrow{d^+} S^-
\xrightarrow{d^-} S^+\xrightarrow{d^+}\cdots
\end{equation}
is exact on
$\EE^0\setminus i(X)$, \cite[\S3.1]{PV}
(exactness is a consequence of a nondegeneracy condition
\cite[\S4.1,~(4.2)]{PV}
on $a$, which is equivalent to Definition
\ref{nondeg}). In \cite{BFM},
a localized Chern character $\ch_X^Y$ is defined in the
rational bivariant
Chow group $A(X\to Y)_\QQ$ for a finite complex of
coherent locally free sheaves on a
scheme $Y$ exact outside a closed subscheme $X$. In \cite{PV},
Section 2.2, this construction is adapted to $2$-periodic complexes
and applied to $S$ to produce a
character $\ch^{\EE^0}_X(S)\in A(X\to \EE^0)_\QQ$,
which allows us to define
\begin{equation}\label{pvctop}
c_{PV}(E^\bullet,b)=\td(E^1)\cdot\ch^{\EE^0}_X(S)\cdot[p]\in A^*(X)_\QQ.
\end{equation}
By the same argument as in Section \ref{sect:derived},
this definition only depends
on $E^\bullet$ and on the homomorphism in
the derived category
$\Sym^{m+1}(E^\bullet)\to \mathcal O_X[-1]$
\linebreak induced by $b$.
Therefore, for any object
$(\pi\colon C\to X;L;f)$ of the stack
$\MSbar$,
we choose $m=r-1$, and we apply
$c_{PV}$ to a complex $E^\bullet$ representing
$R\pi_*L$ and a homomorphism $b$ representing
$\tau\colon \Sym^r(R\pi_*L)\to \mathcal O_X[-1]$.
We obtain the Polishchuk--Vaintrob class
\begin{equation}\label{wtcc-pv}
c_{PV}\in A^*(\MSbar)_\QQ.
\end{equation}

We want to show  $c_{PV}=\cvirt$ in $A^*(\MSbar)_\QQ$.
This amounts to showing
$c_{PV}(E^\bullet,b)=\ctop(E^\bullet,b)$.
First, we show that
Polishchuk and Vaintrob's \linebreak localized
Chern character descends to
a homomorphism from
a $K$-theory group to
the rational bivariant Chow group.

\subsection{Preliminaries on $2$-periodic complexes.}
We need to introduce some notation.
\begin{defn}
Let $Y$ be a scheme and $X$ a closed subscheme.
Write $\mathfrak A$ and $\mathfrak F$ for the categories of
coherent sheaves and coherent locally free sheaves on $Y$.
Denote by ${\calCh}_{\ZZ/2}\mathfrak A$
and  $\calCh_{\ZZ/2}\mathfrak F$
the category of $2$-periodic complexes $W=W^+\oplus W^-$ of coherent
sheaves and coherent locally free sheaves on $Y$.
Write
$\calCh_{\ZZ/2}^X\mathfrak A$ and  $\calCh_{\ZZ/2}^X\mathfrak F$
for the full subcategories of
objects which are
exact on $Y\setminus X$.

We consider the Grothendieck groups
$K_0(\calCh_{\ZZ/2}\mathfrak A)_Y$,
$K_0(\calCh_{\ZZ/2}\mathfrak F)_Y$,\linebreak
$K_0(\calCh_{\ZZ/2}^X\mathfrak A)_Y$,
and  $K_0(\calCh_{\ZZ/2}^X\mathfrak F)_Y$
generated by the objects above
modulo the following relations:
\begin{enumerate}
\item[($i$)]
$[W_1]=[W_2]$ if there exists a quasiisomorphism
$W_1\to W_2 $;
\item[($ii$)]
$[W]=[W_1]+[W_2]$ if there is a sequence
$0\to W_1\to W\to W_2\to 0$ which is exact in
all degrees of $\ZZ/2\ZZ$.
\end{enumerate}
\end{defn}
\begin{rem}\label{rem:PVonK}
Let $Y$ be a smooth scheme
and $X$ a closed subscheme; then,
Polishchuk and Vaintrob's localized Chern character
descends to an homomorphism on
\begin{equation}
\ch_X^Y\colon K_0(\calCh_{\ZZ/2}^X\mathfrak F)_Y
\to A(X\to Y)_\QQ.
\end{equation}
In order to define the homomorphism,
we need to prove that $\ch_X^Y$ is defined
on the objects of $\calCh_{\ZZ/2}^X\mathfrak F$, and is compatible
with $(i)$ quasiisomorphisms and $(ii)$ exact sequences.

In \cite[\S2.2]{PV}, a
technical condition is required for
the existence of \linebreak $\ch_X^Y(W)$:
on $Y\setminus X$, the sheaves
$\im d^+$ and $\im d^-$
have to be subbundles of $W^-$ and $W^+$.
This condition is satisfied, because $Y$ is smooth.
Indeed, for any object $W$ in $\calCh_{\ZZ/2}^X\mathfrak F$,
we can write
resolutions on $Y\setminus X$ of $\im d^+$ (and of $\im d^-$)
of arbitrary length:
$$0\to \im d^+\to W^-_{\mid Y\setminus X}\to
W^+_{\mid Y\setminus X}\to \dots \to W^+_{\mid Y\setminus X}\to
\im d^+\to 0.$$
Take a sequence as above of length $l>\dim(Y)+2$. By
\cite[\S4,~Lem.~9]{BS}, for $Y$ smooth and $W^+, W^-\in \mathfrak F$
we have $\im d^+\in \mathfrak F$.

Finally $\ch_X^Y$ is compatible with
the relations of $K_0(\calCh_{\ZZ/2}^X\mathfrak F)_Y$ because
we have the identity $\ch_X^Y(W)=\ch_X^Y(W_1)+\ch_X^Y(W_2)$
for any exact sequence \linebreak $0\to W_1\to W\to W_2\to 0$ by
\cite[Prop.~2.3,(iv)]{PV}.
Note that this is sufficient by the same argument
of \cite[Exa.~18.1.4]{Fu}: for any
quasiisomorphism $f\colon W_1\to W_2$, we have 
the exact sequence
$0\to W_1\to {\rm Cone}(f)\to  W_2 \to 0$, with
$\ch_X^Y({\rm Cone}(f))=0$
by the exactness of ${\rm Cone}(f)$ on $Y$.
\end{rem}

The following lemma was pointed out
to me by Charles Walter.
\begin{lem}\label{lem:K-2-per}
We have an isomorphism
$K_0(\calCh_{\ZZ/2}\mathfrak A)_Y\to
K_0'(Y)$
\linebreak defined by  $[W]\mapsto [H^+(W)]-[H^-(W)]$.
Furthermore, for $Y$ smooth and quasiprojective,
we also have
$K_0(\calCh_{\ZZ/2}\mathfrak A)_Y\cong
K_0(\calCh_{\ZZ/2}\mathfrak F)_Y$.
\end{lem}
\begin{proof}
The homomorphism
$K_0'(Y)\to K_0(\calCh_{\ZZ/2}\mathfrak A)_Y$
is induced by the functor
sending a sheaf $F$ to the $2$-periodic complex
$F\oplus 0=(\cdots \to F\to 0 \to F\to 0\to \cdots )$.
Composing $[W]\mapsto[H^+(W)]-[H^-(W)]$
after $[F]\mapsto[F\oplus 0]$
we get the identity in
$K_0(\mathfrak A)=K'_0(Y)$.
Conversely, reversing
the order of the composition,
we obtain the homomorphism
$[W]\mapsto[H^+(W)\oplus 0]-[H^-(W)\oplus 0]$.
This functor descends to
the identity in $K_0(\calCh_{\ZZ/2}\mathfrak A)_Y$,
because the sequence in
$\calCh_{\ZZ/2}\mathfrak A$
$$0\longrightarrow \ker d^+\oplus \ker d^-
\longrightarrow W \longrightarrow
\im d^+\oplus \im d^-\longrightarrow 0,$$
is exact ($\ker d^+\oplus \ker d^-$ and
$\im d^+\oplus \im d^-$ are $\ZZ/2\ZZ$-graded
sheaves equipped with the zero differential).

If $X$ is smooth, the $K$-theory groups
$K_0(\calCh_{\ZZ/2}\mathfrak A)_Y$ and
$K_0(\calCh_{\ZZ/2}\mathfrak F)_Y$
are isomorphic, because the
localization of $\calCh_{\ZZ/2}\mathfrak A$
and $\calCh_{\ZZ/2}\mathfrak F$ induces equivalent
triangulated categories.
This follows from the natural inclusion
of $\calCh_{\ZZ/2}\mathfrak F$
into $\calCh_{\ZZ/2}\mathfrak A$ and
the fact that for any object
$W\in \calCh_{\ZZ/2}\mathfrak A$
there exists an exact sequence
$0\to L_n\to L_{n-1}\to \cdots \to L_0\to W\to 0$,
where each $L_i$ is an object of $\calCh_{\ZZ/2}\mathfrak F$.
This is a consequence of the following properties.
\begin{enumerate}
\item \label{tronc} If
$0\to F\to L_p\to L_{p-1}\to \cdots \to L_0\to W\to 0$
is exact with $F, W$ in $\calCh_{\ZZ/2}\mathfrak A$
and $L_i$ in $\calCh_{\ZZ/2}\mathfrak F$; then, for
$p\ge \dim(X)-1$,
$F$ belongs to
$\calCh_{\ZZ/2}\mathfrak F$.
\item \label{freequot}
Any $W$ in
$\calCh_{\ZZ/2}\mathfrak A$ is a quotient of
an object $L$ in $\calCh_{\ZZ/2}\mathfrak F$.
\end{enumerate}
Note that \eqref{tronc} follows immediately from \cite[Lem.~9]{BS}
and \eqref{freequot} can be shown using
the natural functor
\begin{multline}\Hom_{\mathfrak A}(P,W^+)\to\\
\Hom_{\calCh_{\ZZ/2}\mathfrak A}
\left(\left(P\oplus P,\begin{pmatrix}0&\id\\0&0 \end{pmatrix}\right),\left(W^+\oplus W^-,
\begin{pmatrix}0& d^-\\ d^+&0 \end{pmatrix}\right)\right),
\end{multline}
given by $\phi\mapsto \phi \oplus (d_+\circ \phi)$
and the functor
$\phi\mapsto (d_-\circ \phi) \oplus\phi $
defined on
$\Hom_{\mathfrak A}(P,W^-)$.
To see this,
take an object
$W^+\oplus W^-$ in $\calCh_{\ZZ/2}\mathfrak A$,
and, using \cite[Lem.~10]{BS}, take two surjections
$(P^+\to W^+)\in \Hom_{\mathfrak A}(P,W^+)$
and $(P^-\to W^-)\in \Hom_{\mathfrak A}(P,W^-)$.
Then, applying the functors
above, we obtain two morphisms
to $W^+\oplus W^-$.
By construction, the sum is surjective.
\end{proof}

For the rest of the section
we consider a coherent locally free
sheaf $V$ on a smooth scheme $X$, and
we write $\VV$ for the total space of $V$,
$p\colon \VV\to X$ for the projection,
and $i\colon X\to \VV$ for
the zero section.

Note that, for any coherent
$2$-periodic complex $W$ on $\VV$
exact outside $i(X)$,
the sheaves $p_*H^+(W)$ and  $p_*H^-(W)$
are coherent, because
$H^+(W)$ and  $H^-(W)$ are
supported on $i(X)$.
The pushforward $p_*$
induces an exact
functor
from $\calCh_{\ZZ/2}^X\mathfrak A$
to the category of $2$-periodic complexes of coherent sheaves
on $X$.
Therefore, we can
define a homomorphism
\begin{equation}
\fie\colon K_0(\calCh_{\ZZ/2}^X\mathfrak F)_{\VV}\to K_0'(X)
\end{equation}
mapping as $[W]\mapsto
[p_*H^+(W)]-[p_*H^-(W)]$.
\begin{lem}\label{cor:fieiso}
For any vector bundle $\VV$ on a smooth scheme
$X$, the homomorphism $\fie$ is
a bijection.
\end{lem}
\begin{proof}
Note that $K_0(\calCh_{\ZZ/2}^X\mathfrak F)_{\VV}$
is the kernel of the homomorphism
$$K_0(\calCh_{\ZZ/2}\mathfrak F)_{\VV}\to
K_0(\calCh_{\ZZ/2}\mathfrak F)_{\VV\setminus i(X)}.$$
On the other hand the kernel of the homomorphism $K_0'(\VV)\to K_0'(\VV\setminus i(X))$
is identified (via $p_*$) with $K_0'(X)$.
Finally, note that $\fie$ is the restriction
to $K_0(\calCh_{\ZZ/2}^X\mathfrak F)_{\VV}$
of the isomorphism $K_0(\calCh_{\ZZ/2}\mathfrak F)_{\VV}\to
K_0'(\VV)$ from Lemma \ref{lem:K-2-per}.
\end{proof}
\begin{lem}\label{cor:diag}
Let $V$ be a coherent locally free sheaf on $X$.
The diagram
\[\xymatrix{
K_0(\calCh_{\ZZ/2}^X\mathfrak F)_{\VV}
\ar[d]_{\fie}\ar[rr]^{\ch_X^{\VV}} &&
A(X\to \VV)_\QQ\ar[d]^{[p]\cdot\td(V)}\\
K_0'(X) \ar[r]^{\sim}& K_0(X)\ar[r]^{\ch}&A^*(X)_\QQ
} \]
is commutative.
\end{lem}
\begin{proof}
The Koszul complex induced by
the tautological section $\VV\to p^* V$
$$0\to \Lambda^{\rk(V)}(p^*V)^\vee
\xrightarrow{d_{\rk(V)}}
\cdots
\xrightarrow{d_3}
\Lambda^2 (p^*V)^\vee
\xrightarrow{d_2}
(p^*{V})^\vee
\xrightarrow{d_1}
{\mathcal O}_{\VV}
\to 0,$$
is exact off $i(X)$;
this happens in general
for a vector bundle with a nowhere
vanishing section,
see \cite[B.3.4 and A.5]{Fu}.
Denote by $\Lambda (p^*V)^\vee$ the \linebreak $2$-periodic complex
given by summing
the exterior powers $\Lambda^i (p^*V)^\vee$
and the homomorphisms $d_i$ for $i$ even and $i$ odd.
For any coherent locally free sheaf $A$ on $X$,
the $2$-periodic complex
$p^*A\otimes \Lambda(p^*V)$ is
exact off $X$.
We define $\ka\colon
K_0(X)\to K_0(\calCh_{\ZZ/2}^X\mathfrak F)_{\VV}$
as the homomorphism
$[A]\mapsto [p^*A\otimes \Lambda (p^*V)^\vee]$.
Note that $\ka$ commutes with
$\fie$ and $K_0'(X)\to K_0(X)$.
Indeed,
consider the complexes
$A\otimes \Sym^N(V\xrightarrow{\id}V)$,
and denote by $$T_N=(\cdots \to T_N^+\to T_N^-\to T_N^+\to T_N^-\to \cdots)$$
the $2$-periodic complex given by summing over even and odd degrees.
Then, note that,
by the projection formula,
we can decompose the complex
\linebreak $p_*(p^*A\otimes \Lambda (p^*V)^\vee)$
into the sum of the complexes $T_N$.
The complex $T_N$ is exact if $N>0$, so
we have
$$\fie (p^*A\otimes \Lambda p^*V^\vee)=
[H^+(p_*(p^*A\otimes \Lambda p^*V^\vee)]-
[H^-(p_*(p^*A\otimes \Lambda p^*V^\vee)]=[A].$$
Therefore $\ka$ is an isomorphism.

By \cite[Prop.~2.3,~(iv)]{PV},
we have
$$\ch^{\VV}_X(p^*A\otimes \Lambda
(p^*V)^\vee)
\cdot[p]\cdot \td(V)=\ch(A)$$
(indeed, the formula above can be rewritten
in terms of the standard localized Chern character from
\cite{Fu} and follows from
\cite[Prop.~3.4, Ch.~1]{BFM}).
The diagram commutes, because we have
$$\ch^{\VV}_X(\ka[A])\cdot[p]\cdot \td(V)=
\ch^{\VV}_X(p^*A\otimes \Lambda
(p^*V)^\vee)\cdot[p]\cdot \td(V)=\ch(A).$$
\end{proof}

\subsection{The identity.}
We show the identity between $c_{PV}$ and $\cvirt$.
For any object
$(\pi\colon C\to X;L;f)$ of the stack
$\MSbar$,
we have a complex $E^\bullet$ and a homomorphism
$b\colon \Sym ^{m}E^0\otimes E^1\to \mathcal O_X$,
whose dual is closed and nondegenerate for $(E^\bullet)^\vee$.
We recall that
$c_{PV}(F_\bullet,a)=\td (E^1)\cdot\ch^{\EE^0}_X(S)\cdot[p]$,
where $S$ is the $2$-periodic complex at
\eqref{P-cx}.
\begin{thm}\label{thm:compPV}
We have \begin{align}
&\  &&\fie[S]=\ke((E^\bullet)^\vee,b^\vee)
&&\text{in $K_0'(X)$ and}\label{K-compPV}\\
&\ &&c_{PV}(E^\bullet,b)=\ctop(E^\bullet,b)
&&\text{in $A^*(X)_\QQ$}\label{A-compPV}.
\end{align}
Therefore
$c_{PV}=\cvirt$ in $A(\MSbar)_\QQ$.
\end{thm}
\begin{proof}
Define $L^+\to L^-$ and $L^-\to L^+$ by summing,
for $k$ even and $k$ odd, all the homomorphisms $\wt d:L^{h,k}\to L^{h+1,k-1}$ and
$\wt a\colon L^{h,k}\to L^{h+m,k+1}$.
Now, $L=L^+\oplus L^-$ is a $2$-periodic complex of
quasicoherent sheaves by Lemma \ref{doublecx}.
It is easy to see that the $2$-periodic complex
$p_*S$ is equal to $L$ (with \linebreak $p_* \widetilde \delta=\widetilde d$
and $p_*\widetilde \fie=\widetilde a$).
The functor $p_*$ is exact on complexes which are
exact off $i(X)$; therefore,
since the complex $S$ in \eqref{P-cx}
is exact off $i(X)$ by \cite[S3.1]{PV},
we have
$p_*\left({ H}^+(S)\right)= H^+(p_*S)={ H}^+(L)$
and $p_*\left({H}^-(S)\right)=H^+(p_*S)={H}^-(L)$.
This implies \eqref{K-compPV}.
Finally, Lemma \ref{cor:diag} implies \eqref{A-compPV}:
\begin{multline*}
c_{PV}(E^\bullet,b)=\td (E^1)\ch^{\EE^0}_X(S)\cdot[p]=
\td(E^1)\cdot \td (E^0)^{-1}\cdot \ch(\fie[S])\\=
\frac{\td(E^1)}{\td(E^0)}
\ch(\ke((E^\bullet)^\vee,b^\bullet))=
\ctop(E^\bullet,b).
\end{multline*}
\end{proof}
\begin{rem}\label{propPV}
The identity of Theorem \ref{thm:compPV}
guarantees that the class is concentrated in
degree $-\chi(C,L)=n_1-n_0$.
Furthermore, the axioms of
cohomological field theory stated
by Jarvis, Kimura, and Vaintrob in
\cite{JKV} are proven by Polishchuk and Vaintrob in
\cite{PV} and \cite{Po}.
The proofs are given on the level of $K$-theory of complexes
on a vector bundle that are strictly exact off the zero section.
Note
that the equivalence of
triangulated categories shown in  Lemma 5.3.4
implies that working in the $K$-theory of the base
scheme $X$ is equivalent to working in the $K$-theory
of complexes on a vector bundle
$\VV\to X$ that are strictly exact off the zero section.
This indicates that the arguments used by Polishchuk
and Vaintrob can be restated in $K_0(X)$,
without passing through $\VV$.
In the cases that we have checked, however,
this passage to $X$
does not appear to simplify
significantly the proofs of the cohomological field theory axioms.
\end{rem}

\section{Computations}
\label{sect:comp}
\subsection{Compatibility with Witten's definition.}
Assume that for every point $x\in X$ we have $H^0(C_x,L_x)=0$.
Then, $R^0\pi_*L=0$ and $R^1\pi_*L$ is a bundle.
By Definition \ref{wtcc-chow},
we have $\cvirt=\ctop(R^\bullet\pi_*L,0)$.
By Definition \ref{defnctop}, we have
$$\ctop(R^\bullet\pi_*L,0)=\ch(\ke((R^\bullet\pi_*L)^\vee, 0))\cdot
{\td(R^1\pi_*L}).$$
So, by Proposition \ref{pro:oneiszero}, we have
$$\ke((R^\bullet\pi_*L)^\vee, 0)=\lambda_{-1}((R^1\pi_*L)^\vee).$$
Finally, by \eqref{ke-ctop},
we have
$$\cvirt=\ctop(R^\bullet\pi_*L,0)=
\lambda_{-1}((R^1\pi_*L)^\vee)\cdot
{\td(R^1\pi_*L})=\ctop(R^1\pi_*L),$$
which agrees with Witten's definition
\eqref{wittensdef}.

\subsection{The case when $R^0\pi_* L$ and $R^1\pi_* L$
are vector bundles.}
Let $(\pi\colon C\to X; L; f)$ be an object of $\MSbar$
over which $R^ip_*L$ is a vector bundle
of rank $h^i=h^i(C,L)$.
Denote by $a$ the form induced by Serre duality
$$b\colon \Sym^{r-1} (R^0p_*L)\otimes (R^1p_*L)\to \mathcal O_X.$$
By the same argument as in Section \ref{sect:wtcc}, the form
$b^\vee$ is closed and nondegenerate.
By Theorem \ref{bound} we have, for $t_0=r-1+\binom{h^1-1+r-1}{r-1}$,
$$\ke (Rp_*L^\vee,b^\vee)=\sum_{h\le (r-1)k+t_0}
(-1)^k [\Sym^h (R^0p_*L)^\vee][\Lambda^k (R^1p_*L)^\vee].$$
Therefore, via Definition \ref{defnctop}, we
can compute explicitly the
pullback of $c_W$ under the morphism $X\to \MSbar$.

\subsection{The case $r = 2$, theta characteristics.}
For $r = 2$ and $\kb=\pmb 0$, the class $\cvirt$ is
a cycle of
codimension $0$, Remark \ref{propPV}.
So the class $\cvirt$ is determined by a number.
We claim that this number is $1$ on the connected component
compactifying even theta characteristics ($h^0(C,L)$ even)
and $-1$ on the component compactifying the odd ones
($h^0(C,L)$ odd).
Indeed, since we only need to determine the
multiplicity of $\cvirt$
we can consider $L\to C\xrightarrow {\pi} X$ where $X$ is a point,
$C$ is a smooth curve, and $L$ is a line bundle
satisfying $L^{\otimes 2}\cong \omega_C$.
The complex $H^\bullet(C,L)$ with zero
differential represents $R\pi_* L$.
The nondegenerate form is the perfect pairing
$H^0(C,L )\otimes H^1(C, L)\to \CC$.
We write $H^\vee_i=(H^i(C,L))^\vee$ and $h^i=h^i(C,L)$
and describe $L^{\bullet,\bullet}$
as follows.
\[ \xymatrix@R=0.1cm{
0& 0& \cdots&\\
\Lambda^{h^1}H^\vee_1 &
    H^\vee_0\otimes\Lambda^{h^1}H^\vee_1 &
      \cdots \\
\Lambda^{h^1-1}H^\vee_1  \ar[ru]_{\wt a}&\cdots
& \Lambda^2 H^\vee_0\otimes H^\vee_1 \\
%\Lambda^{2}H^\vee_1 &
%    H^\vee_0\otimes\Lambda^{2}H^\vee_1 &
%      \cdots &
%            \Sym^{t}{H^\vee_0}\otimes\Lambda^{2}{H_1^\vee} &
%                \cdots \\
\cdots &
    \ H^\vee_0\otimes H^\vee_1 \ \ar[ru]_{\wt a} &
      \Lambda^2 H^\vee_0\otimes H^\vee_1 & \cdots \\
\quad \ \ \CC\quad  \ \  \ar[ru]_{\wt a}&
    \quad \  H^\vee_0 \quad \ \ar[ru]_{\wt a} &
      \cdots &
            \Sym^{t}{H^\vee_0}&
                \cdots \\
} \]
Note that the differential $\wt d$ is zero and the differential
$\wt a$ is always exact, except on $\Lambda^{h^1}H^\vee_1$.
The Witten top Chern class in $K$-theory is
equal to $(-1)^{h^1}[\Lambda^{h^1}H^\vee_1]=
(-1)^{h^1}[\CC]$. The cohomology
class is $\cvirt=(-1)^{h^1}$.
\subsection{Genus one.}
We consider the case $g=1$, $n=1$, and $\kb=(0)$, the
moduli stack compactifying elliptic $r$-spin curves.
Indeed, again the class $\cvirt$ is a cycle of
codimension $0$, Remark \ref{propPV}.
We can consider $L\to C\xrightarrow {\pi} X$ where $X$ is a point,
$C$ is a smooth curve of genus $1$, and $L$ is a line bundle
satisfying $L^{\otimes r}\cong \omega_C\cong \mathcal O_C$.
The complex $H^\bullet(C,L)$ with zero
differential represents $R\pi_* L$.
Note that, if $L$ is not trivial, we have $H^i(C,L)=0$ and, therefore,
$c_W=1_X$.

We assume that $L$ is trivial; therefore, we have $h^i(C,L)=1$.
We write
$H^i=H^i(C,L)$ and $H_i^\vee=(H^i(C,L))^\vee$.
The closed and nondegenerate
form is given by the isomorphism $\Sym^{r-1}H^0
=(H^0)^{\otimes r-1}
\cong H^0(L^{\otimes r-1})\cong H^\vee_1$.
We describe $L^{\bullet,\bullet}$
as follows.
\[ \xymatrix@R=0.1cm{
0& 0& \cdots&\\
H^\vee_1 &
    H^\vee_0\otimes H^\vee_1 &
        \cdots&
            (H^\vee_0)^{\otimes r-2}\otimes H^\vee_1 &
                (H^\vee_0)^{\otimes r}\otimes H^\vee_1 \\
\CC&
    H^\vee_0 &
        \cdots&
            (H^\vee_0)^{\otimes r-2} &
                (H^\vee_0)^{\otimes r} \\
} \]
The differential $\wt d$ is zero, and the differential
$$\wt a\colon (H^\vee_0)^{\otimes i}\to
(H^\vee_0)^{\otimes i+r-1}\otimes H^\vee_1$$
is always exact, except on
$H^\vee_1$,  $H^\vee_0\otimes H^\vee_1$, \dots, and
$(H^\vee_0)^{\otimes r-2}\otimes H^\vee_1$.
Therefore,
the
Witten top Chern class in the $K$-theory
of a closed point 
is equal to $-\sum_{i=1}^{r-2}[(H^\vee_0)^{\otimes i}\otimes H^\vee_1]=
(-1)[\CC^{\oplus r-1}]$.
Finally, via Chern character, we get $\cvirt=-(r-1)1_X$.
This is consistent with Witten's predictions
and can be deducted also from
the cohomological field theory axioms \cite{JKV}.

\nocite{*}
\bibliographystyle{amsalpha}
\bibliography{bibchiodo}

\end{document}